\documentclass{elsart}

\usepackage{amssymb}
\usepackage[english,francais]{babel}

\newtheorem{theorem}{Theorem}[section]

\newtheorem{definition}[theorem]{Definition\rm}

\newcommand{\Vol}{{\mathrm{Vol}}}
\newcommand{\ehr}{{\mathrm{ehr}}}
\newcommand{\Z}{{\mathbb Z}}
\def\og{\leavevmode\raise.3ex\hbox{$\scriptscriptstyle\langle\!\langle$~}}
\def\fg{\leavevmode\raise.3ex\hbox{~$\!\scriptscriptstyle\,\rangle\!\rangle$}}

\journal{the Acad\'emie des sciences}
\begin{document}
\centerline{Combinatorics}
\begin{frontmatter}


\title{Une relation entre nombre de points entiers, volumes des faces et degr\'e du
discriminant des polytopes entiers non singuliers}
 \author{Alicia Dickenstein,}
 \address{Departamento de Matem\'atica,
FCEN,
Universidad de Buenos Aires and IMAS, CONICET,
Ciudad Universitaria,  Pab I,
(C1428EGA) Buenos Aires, Argentina}
 \ead{alidick@dm.uba.ar}
\thanks{AD is partially supported by UBACYT 20020100100242, CONICET
 PIP 112-200801-00483 and ANPCyT 2008-0902, Argentina}
 \author{Benjamin Nill,}
\address{Case Western Reserve University, Department of Mathematics, 10900 Euclid Avenue, Cleveland, OH 44106, USA}
 \ead{benjamin.nill@case.edu}
\thanks{BN is supported by the US National Science Foundation (DMS 1203162)}
 \author{Mich\`ele Vergne}
\address{Institut de Math\'ematiques de Jussieu, 175 rue du Chevaleret, 75013 Paris}
 \ead{vergne@math.jussieu.fr}

\selectlanguage{english}
\title{A relation between number of integral points,
volumes of faces and degree of the discriminant of smooth lattice polytopes}

\begin{center}
{\small Received *****; accepted  +++++\\}
\end{center}

\begin{abstract}
\selectlanguage{english}
We present a formula for the degree of the discriminant
of a smooth projective toric variety associated to a lattice
polytope $P$, in terms of the number
of integral points in the interior of
dilates of faces of dimension
greater or equal than  $\lceil \frac {\dim P} 2 \rceil$.
 {\it To cite this article: A. Dickenstein, B. Nill, M.
 Vergne, C. R. Acad. Sci. Paris, Ser. I  (2012).}

\vskip 0.5\baselineskip

\selectlanguage{francais}
\noindent{\bf R\'esum\'e} \vskip 0.5\baselineskip \noindent

Nous donnons une formule pour le degr\'e du discriminant
d'une  vari\'et\'e torique projective non singuli\`ere associ\'ee \`a
un polytope entier $P$, en terme du nombre de points
entiers des int\'erieurs de
dilatations de faces de dimension sup\'erieure
ou \'egale \`a $\lceil \frac {\dim P} 2 \rceil$.
Pour citer cet article~:  A. Dickenstein, B. Nill, M. Vergne,
C. R. Acad. Sci. Paris, Ser. I  (2012).
\end{abstract}

\end{frontmatter}

\section{Introduction}
The relation between the volume and the number of integral points of a lattice polytope
$P \subset {\mathbb R}^n$ has a long history. Here, a lattice polytope is a polytope
whose vertices are integral, i.e. in $\Z^n$. The  Ehrhart polynomial \cite{E}
$\ehr_P(t)$ is  the  polynomial  of degree $\dim (P)$
 in one variable $t$ such that
the number of integral points in $tP$ is equal to $\ehr_P(t)$
when $t$ is a non negative integer.
On the other hand, by Ehrhart-Macdonald reciprocity, for $t$ a
positive integer, $\ehr_P(-t)$ equals
$(-1)^{\dim(P)}$ times the number of  integral points
in the relative interior of $tP$. The leading coefficient
of these polynomials equals $1/{\dim(P)!}$ times the lattice
volume $\Vol_{\mathbb Z}(P)$ of $P$.
(The volume $\Vol_{\mathbb Z}$ is normalized so that the volume
of the fundamental parallelepiped  is $\dim(P)!$.)

An $n$-dimensional lattice polytope $P$ with $N+1$ integral points
defines an embedded projective toric variety $M_P \subset {\mathbb P}^N$.
The polytope $P$ is called smooth (or Delzant)
when $M_P$ is nonsingular, which implies that $P$ is simple and the primitive
vectors on the $n$ edges emanating from each vertex form
a basis of ${\mathbb Z}^n$.
The dual projective variety $M_P^\vee \subset ({\mathbb P}^N)^\vee$
consisting of the closure of the locus of those hyperplanes that intersect $M_P$ non-transversally is
generically a hypersurface. When this is the case, its degree equals
\begin{equation}\label{eq:vol}
c(P)=\sum_{p=0}^n (-1)^{n-p}(p+1)\sum_{F\in \mathcal F_p(P)} \Vol_{\mathbb Z}(F),
\end{equation}
where $\mathcal F_p(P)$ denotes the subset of $p$-dimensional faces of $P$
\cite[Th.~28, Ch.~9.2]{GKZ}. In fact $c(P)$ equals the top Chern class of the first jet bundle of
the embedding, and the lattice volumes in (\ref{eq:vol}) occur as combinatorial
translations of intersection products of line bundles on $M_P$.

In this note, our main result, Theorem~\ref{cor:degree}, gives a new representation of $c(P)$
for any simple $n$-dimensional lattice polytope $P$
in terms of the number of integral points in the interior of 
dilates of faces of dimension
greater or equal than  $\lceil \frac n 2 \rceil$.

The search for this representation was motivated by the question by Batyrev and Nill raised in \cite{BN},
whether a lattice polytope with sufficiently large dilates without interior lattice points necessarily has a
Cayley structure. While this question was answered affirmatively in \cite{HNP}, it is in general still open what
`sufficiently large' precisely means. In the smooth case, this has recently been clarified in \cite[Th.~2.1 (i)]{DN}:
if the $\frac n 2 + 1$ dilate of a smooth polytope $P$ does not have interior lattice points, then $c(P) = 0$, so
$M_P^\vee$ is not a hypersurface ($M_P$ is {\em{dual defective}}), and
$P$ admits a (strict) Cayley structure (cf. \cite{DR}).
The Ehrhart-theoretic proof in \cite{DN} relied heavily on non-trivial binomial identities and lacked any general insight in why $\frac n 2 + 1$ works.
Moreover, the method of proof didn't employ Ehrhart reciprocity and didn't allow to deal with lattice polytopes with interior lattice points.

In Section~\ref{section1}, we provide a general result, Theorem~\ref{Nill}, on involutions of `Dehn-Sommerville type'.
The proof of this formal statement is elementary and short. As our main application we
deduce from Ehrhart reciprocity our new formula for $c(P)$ valid for any simple lattice polytope (not necessarily smooth).
This shows that
the surprising fact that $c(P)$ depends only on the lattice points of dilates of faces of high dimension
holds actually for {\em any} simple polytope.

In Section~\ref{sec:applications}, we explain how applying Theorem~\ref{Nill}
yields the desired direct and conceptual proof of \cite[Th.~2.1 (i)]{DN}.


We show in the last section an application of Theorem~\ref{Nill} in the realm of symmetric functions. We thank Alain Lascoux and Jean-Yves Thibon
for their explanations about the ubiquitous occurrence of the transformation $S$ in Definition~\ref{defS}, see for example   \cite{thibon}.

Finally, let us remark that the method of proof of Theorem~\ref{Nill} can also be used to
show an expression of the volume of $P$ in terms
of the number of integral points and integral boundary points
on $\le \lceil \frac n 2 \rceil$-dilates of $P$, similar to \cite[Th.~2.5]{Kas}.


%

\section{An  involution} \label{section1}
\begin{definition} \label{defS}
Define the  transformation $S:{\mathbb C}^{n+1}\to {\mathbb C}^{n+1}$
by $$S[E_0,\ldots, E_n]=[F_0,\ldots, F_n]$$ with
$$F_p=\sum_{j=0}^p (-1)^j {n-j \choose n-p}E_j.$$
\end{definition}

Then $S^2=1$, and we have the identity
\begin{equation}\label{magique}
\sum_{p=0}^n(-z)^p(z+1)^{(n-p)}E_p=\sum_{p=0}^n  z^pF_p.
\end{equation}
Later on, the elements $E_i$ will be  themselves functions of other variables.

The definition of $S$  is motivated by the following examples.

\medskip

{\bf Example 1.} Let $[f_0, \ldots, f_n]$ be the $f$-vector
of a simple $n$-dimensional polytope.
Then  $S([f_0, \ldots, f_n])=$ $[f_0, \ldots, f_n]$ by the Dehn-Somme\-rville
equa\-tions, see e.g., \cite[Th.~5.1]{BeRo07}.
On the other hand, in the dual situation when $P$ is a simplicial $n$-dimensional polytope
(i.e., all the facets of $P$ are simplices), $S$ is (up to a sign) precisely
the transformation between the $f$- and $h$-vectors, see \cite[8.3]{Zie}.

\medskip

{\bf Example 2.} Let $P\subset {\mathbb R}^n$ be a polytope with integral vertices.
 Denote by $\mathcal F_k(P)$
 the set of   $k$-dimensional faces of $P$ and define
 $E_k^P(t)=\sum_{F\in \mathcal F_k(P)} \ehr_F(t)$. We obtain an element
 ${\bf E}^P=$ $[E_0^P(t),E_1^P(t),\ldots, E_n^P(t)]$  of ${\mathbb C}^{n+1}$
 depending of $t$. We write $S({\bf E}^P)=[F^P_1(t),F^P_2(t),\ldots, F^P_n(t)]$.

 Assume $P$ is simple. In this case, any face $F$ of $P$ of
 dimension $j \le p$ is contained in precisely ${n-j \choose n-p}$
faces $G$ of $P$ of dimension $p$.
Using the inclusion-exclusion formula, we get, for every $p=0, \dots, n$, the extended Dehn-Sommer\-ville equations
\[ E_p^P(-t) = \sum_{j=0}^p (-1)^{j} {n-j \choose n-p} E^P_j(t) =  F^P_p(t).\]
Thus, for any non negative integer $t$, $(-1)^p F^P_p(t)$
equals the sum of the number of integral
points in the relative interior of the $p$-dimensional faces of $tP$.
This example is Theorem~5.3 in \cite{BeRo07}.

\medskip

{\bf Example 3.}
The elementary symmetric functions $\sigma_i({\bf x})$, ${\bf x} =(x_1, \dots, x_n)$,
are defined by
$\prod_{a=1}^n(x_a z+1)=\sum_{i=0}^n \sigma_i({\bf x}) z^i$. If we start with
the sequence $E=[\sigma_0({\bf x}),\ldots,\sigma_n({\bf x})]$, then $S(E)$ gives the sequence of
elementary functions on $(1-x_1,\ldots,1-x_n)$. We clearly see in this case
that $S$ is an involution.
More generally, let $v,w$ be meromorphic functions of one variable with $v(x)+w(-x)=1$,
write $\prod_{a=1}^n(v(x_a) z+1)=\sum_{i=0}^n V_i({\bf x}) z^i$,
 $\prod_{a=1}^n(w(x_a) z+1)=\sum_{i=0}^n W_i({\bf x}) z^i$.
 If we start with
the sequence $E=[V_0({\bf x}),\ldots,V_n({\bf x})]$, then
$S(E)=[W_0({-\bf x}),\ldots, W_n({-\bf x})]$ .

%

\medskip

Let $\mathcal P_n$ be the space of families $[E_0(t), \ldots, E_n(t)]$ of $n+1$ polynomials
 where $E_j(t)$ is a polynomial in $t$ of degree less or equal to $j$.
The transformation $S$ induces a transformation of $\mathcal P_n$.
We write each element $
E_j(t)= \frac{t^j}{ j!} v_j+\mbox{lower  terms}$.
 Define
 \begin{equation}\label{eq:c}
 c({\bf E})=\sum_{p=0}^n (-1)^{n-p} (p+1) v_p.
 \end{equation}

\begin{theorem}\label{Nill}
Let ${\bf E}=[E_0,E_1,\ldots, E_n]\in \mathcal P_n$ and ${\bf SE}=[F_0,F_1,\ldots,F_n]$.
For $n$ odd, and $m=(n+1)/2$,
then
\begin{equation}\label{eq:odd}
c({\bf E})=\sum_{p=m}^n\sum_{i=1}^{p+1-m}(-1)^{p+m-i}
 {p+1 \choose m+i} i (E_p(-i)+ F_p(i)).
\end{equation}
For $n$ even, and $m=n/2$, then
{\small
\begin{equation}\label{eq:even}
\! \! \! \! \! \! \! \! c({\bf E})=\sum_{p=m}^n\sum_{i=1}^{p+1-m}(-1)^{p+1+m-i}
 \left({p+1 \choose m+i}- {p+1 \choose m+i+1}\right)
\frac{i}{2} (E_p(-i)+F_p(i)).
\end{equation}
}
\end{theorem}

Let us give the proof of this identity for $n$ odd, the case $n$ even being similar.

Denote by $\tau$ be the translation operator $(\tau h)(t)=h(t+1)$.
Define $e_p(t)=tE_p(t)$, a polynomial function of degree $p+1$.
Then  $(\tau-1)^{p+1} e_p$ is just
the constant function $(p+1)v_{p}$, as can be checked on binomials.
 Hence, $(\tau^{-1}-1)^{p+1} e_p = (-1)^{p+1} (p+1) v_p$.
Since translating by $\tau^m$ doesn't change a constant function, we obtain that
\[c({\bf E}) =   \sum_{p=0}^n (-1)^{n+1} \tau^{m}(\tau^{-1}-1)^{p+1}e_p(0).\]
Since $\tau^{m}(\tau^{-1}-1)^{p+1}$ equals
\[ \sum_{j=0}^{n+1} {p+1 \choose j} \tau^{-j+m} (-1)^{p+1-j} =
\sum_{i=-m}^m {p+1 \choose m-i} \tau^i (-1)^{p+1+i-m},\]
we get
{\small
\begin{equation}\label{true}
 c({\bf E})=\sum_{p=0}^n\sum_{i=-m}^m (-1)^{p+1+i-m} {p+1 \choose m-i} i\, E_p(i).
 \end{equation}
}
Now we write the right hand side of~(\ref{eq:odd}) as
$$RHS:=\sum_{p=m}^n\sum_{i=0}^{p+1-m} (-1)^{p}\, i\,{\rm coeff}((1-z)^{p+1},z^{m+i})
 \,(E_p(-i)+F_p(i)).$$
If $i>p+1-m$,   or $p<m$, the number   $i({\rm coeff}(1-z)^{p+1},z^{m+i})$
is equal to $0$. Thus, since $p+1-m\le n+1-m \le m$,
$$RHS = -\sum_{i=0}^m i \sum_{p=0}^{n} {\rm coeff}((z-1)^{p+1},z^{m+i})
\,(E_p(-i)+F_p(i)).$$
By Relation (\ref{magique}) we get $\sum_p (z-1)^p F_p(i)=\sum_p (1-z)^p\,
z^{(n-p)}E_p(i)$; hence,
$\sum_p (z-1)^{p+1} F_p(i)=-\sum_p (1-z)^{p+1}\, z^{(n-p)}E_p(i)$. We deduce that RHS equals
{\small
\[\! \!  \! \! \! \! \!  \! \! \!\sum_{p=0}^n \sum_{i=0}^m \left({\rm coeff}((z-1)^{p+1},z^{m-(-i)}) \, (-i)\,
E_p(-i)+ {\rm coeff}((1-z)^{p+1} ,z^{m+i-(n-p)}) \, i \,E_p(i)\right).\]}
As {\small ${p+1 \choose m+i-n+p} = {p+1 \choose m-i}$},
this is equal to the expression of $c({\bf E})$ given in Equation~(\ref{true}).

\section{An application to lattice polytopes}\label{sec:applications}
We return to Example~2.
Let us denote by $I_p(i)$ the number of integral points in the relative interior of
$i$-th multiples of $p$-dimensional faces of $P$.

\begin{theorem}\label{cor:degree}
Let $P$ be an $n$-dimensional simple lattice polytope.
Let $$c(P)=\sum_{p=0}^n (-1)^{n-p}(p+1)\sum_{F\in \mathcal F_p(P)} \Vol_{\mathbb Z}(F).$$
For $n$ odd, and $m=(n+1)/2$,
then $$c(P)=\sum_{p=m}^n\sum_{i=1}^{p+1-m}(-1)^{m-i}   {p+1 \choose m+i} 2 \,i \, I_p(i).$$

For $n$ even, and $m=n/2$,
then
 $$c(P)=\sum_{p=m}^n\sum_{i=1}^{p+1-m}(-1)^{m+1-i}  \left({p+1 \choose m+i}-{p+1 \choose m+i+1}\right)
\, i\, I_p(i).$$
\end{theorem}

This result follows from Theorem~\ref{Nill}. Indeed, let ${\bf E}^P$ be the sequence of polynomials described in Example~2.
Then, the  coefficient $v_p$
equals the normalized volume of the skeleton of $p$-dimensional faces of $P$,
and  $F^P_p(i) = E^P_p(-i) = (-1)^p I_p(i)$.

In particular, we get the following alternative proof of \cite[Th.~2.1 (i)]{DN}.
Assume $I_n(i) = \ehr_P(-i)=0$ for any positive integer $i \le \frac n 2 +1$.
For simplicity, take $n$ odd. So, the polytope $iP$ has no integral interior points for
$i=1,2,\ldots, m= \frac {n+1}2$.
The monotonicity theorem of Stanley \cite{S} implies that any face of codimension $k$  of $P$ has no
integral interior points for $i=1, \dots, m-k$.
Thus we obtain from~(\ref{eq:odd}) that  $c(P) = c({\bf E})=0$. If $P$ is smooth, then this implies that $M_P$ is dual defective.

\section{An identity of symmetric functions}
Let
\[V(s,{\bf x})(t)=\frac{e^{ts}}{\prod_{a=1}^n  (-t x_a)}, \quad \quad
B(s,{\bf x})(t)= \frac{e^{ts}}{\prod_{a=1}^n (1-e^{tx_a})}\]
be meromorphic functions of $t$ depending of the $ (n+1)$ variables $s$ and ${\bf x}$.
The constant term $CTV(s,{\bf x})$ of the Laurent series (in $t$) of
 $V(s,{\bf x})(t)$
 is $(-1)^n\frac{s^n}{n!x_1x_2\cdots x_n}.$
The constant term $CTB(s,{\bf x})$ of the Laurent series of
 $B(s,{\bf x})(t)$
 is a meromorphic function of $(s,x_1,x_2,\ldots,x_n)$, symmetric in the $x_i$.

Let $P\subset {\mathbb R}^n$ be an $n$-dimensional smooth polytope, and  let $\mathcal V(P)$ be the set of vertices of $P$.
For $v$ a vertex,  let $g_a,a=1,\ldots,n,$ be the primitive generators on the edges of $P$ starting at $v$.
 If $\xi$ is generic in the dual space to ${\mathbb R}^n$, we can specialize the variables $s,x_a$ to $\langle v,\xi\rangle $, $\langle g_a,\xi\rangle $ in
 $V(s,{\bf x})(t)$  and  $B(s,{\bf x})(t)$ and we obtain  meromorphic functions $V_v(\xi,t), B_v(\xi,t)$ depending of the vertex $v$.
Then the sum over $v\in \mathcal V(P)$ of the meromorphic function $V_v(\xi,t)$ or $B_v(\xi,t)$ is actually regular at $t=0$, and
$n!\sum_v V_v(\xi,t)_{(t=0)}$ is the normalized volume of  $P$ while
$\sum_v B_v(\xi,t)_{(t=0)}$ is the number of integral points in $P$ by Brion's formulae \cite{B}.

 Let $1\leq p\leq n$,
 and let
 \[V_p(s, {\bf x})(t)=\sum_{J, |J|=p}\frac{e^{ts}}{\prod_{a\in J}(-tx_a)}, \quad \quad
 B_p(s,{\bf x})(t)=\sum_{J, |J|=p}\frac{e^{ts}}{\prod_{a\in J}(1-e^{tx_a})},\]
 where $J$ runs over subsets of cardinal $p$ of $\{1,2,\ldots,n\}$. Similarly, the specialization
 $s,x_a$ to $\langle v,\xi\rangle $, $\langle g_a,\xi\rangle $ gives us  meromorphic functions of $(t,\xi)$,
 and the sum over the vertices $v$ of the polytope $P$  is regular at $t=0$.
 If we dilate $P$ by $i$, the vertices are changed in $iv$, while the generators $g_a$ stay the same.
 Then the identities  of Corollary \ref{cor:degree} imply in particular that for $n$ odd,
$$\sum_{v\in \mathcal S(P)} \left(\sum_{p=0}^n (-1)^{n-p} (p+1)! V_{v,p}(\xi,t)\right)_{(t=0)}$$ $$=\sum_{v\in \mathcal S(P)}
\left(\sum_{p=m}^n(-1)^p\sum_{i=1}^{p+1-m}(-1)^{m-i}   {p+1 \choose m+i} 2 \,i \, B_{iv,p}(-\xi,t)\right)_{(t=0)}.$$
A similar identity holds for $n$ odd.
   Actually, this identity holds before summing over the vertices $v$ and before specializing,
   as an identity  for the symmetric functions in $x_a$
     obtained as  the constant term
    $CTB_p(s, {\bf x})$ of the Laurent series in $t$  of
   $B_p(s, {\bf x})(t)$.

\begin{theorem}
For $n$ odd,   and $m=(n+1)/2$,
$$\sum_{p=0}^n (-1)^{n-p}(p+1)! CTV_p(s,{\bf x})= $$
$$\sum_{p=m}^n(-1)^p\sum_{i=1}^{p+1-m}(-1)^{m-i}   {p+1 \choose m+i} 2 \,i \, CTB_p(-i s,{\bf x}).$$

For $n$ even, and  $m=n/2$,
$$\sum_{p=0}^n (-1)^{n-p} (p+1)! CTV_p(s,{\bf x})=$$
$$\sum_{p=m}^n(-1)^p\sum_{i=1}^{p+1-m}(-1)^{m+1-i}  \left({p+1 \choose m+i}-{p+1 \choose m+i+1}\right) \,i \, CTB_p(-i s,{\bf x}).$$
\end{theorem}
We prove this again as a consequence of Theorem~\ref{Nill} using Example 3, and the identity $\frac{1}{1-\exp(x)}+\frac{1}{1-\exp(-x)}=1.$

\medskip

\noindent {\bf Acknowledgements:} We thank the Institut Mittag Leffler,
where this work was initiated. MV gratefully acknowledges support from the AXA Mittag-Leffler
Fellowship Project, funded by the AXA Research Fund.

\end{document}